\documentclass[sn-basic]{sn-jnl}

\usepackage{graphicx}%
\usepackage{multirow}%
\usepackage{amsmath,amssymb,amsfonts}%
\usepackage{amsthm}%
\usepackage{mathrsfs}%
\usepackage[title]{appendix}%
\usepackage{xcolor}%
\usepackage{textcomp}%
\usepackage{manyfoot}%
\usepackage{booktabs}%
\usepackage{algorithm}%
\usepackage{algorithmicx}%
\usepackage{algpseudocode}%
\usepackage{listings}%
\usepackage{subcaption}

\theoremstyle{thmstyleone}%
\theoremstyle{thmstyletwo}%

\theoremstyle{thmstylethree}%

\begin{document}

\title[Article Title]{Data-driven topology design with persistent homology for enhancing population diversity}


\author[1]{\fnm{Taisei} \sur{Kii}}\email{kii@syd.mech.eng.osaka-u.ac.jp}

\author*[1]{\fnm{Kentaro} \sur{Yaji}}\email{yaji@mech.eng.osaka-u.ac.jp}

\author[2]{\fnm{Hiroshi} \sur{Teramoto}}\email{teramoto@kansai-u.ac.jp}

\author[1]{\fnm{Kikuo} \sur{Fujita}}\email{fujita@mech.eng.osaka-u.ac.jp}

\affil[1]{\orgdiv{Department of Mechanical Engineering}, \orgname{Osaka University},\\ \orgaddress{\street{2-1 Yamadaoka}, \city{Suita, Osaka}, \postcode{565-0871}, \country{Japan}}}

\affil[2]{\orgdiv{Department of Mathematics}, \orgname{Kansai University},\\ \orgaddress{\street{3-3-35 Yamate-cho}, \city{Suita, Osaka}, \postcode{564-8680}, \country{Japan}}}


\abstract{
    This paper proposes a selection strategy for enhancing population diversity in data-driven topology design (DDTD), a topology optimization framework based on evolutionary algorithms (EAs) using a deep generative model.
    While population diversity is essential for global search with EAs, conventional selection operators that preserve diverse solutions based on objective values may still lead to a loss of population diversity in topology optimization problems due to the high dimensionality of design variable space and strong nonlinearity of evaluation functions.
    Motivated by the idea that topology is what characterizes the inherent diversity among material distributions, we employ a topological data analysis method called persistent homology.
    As a specific operation, a Wasserstein distance sorting between persistence diagrams is introduced into a selection algorithm to maintain the intrinsic population diversity. 
    We apply the proposed selection operation incorporated into DDTD to a stress-based topology optimization problem as a numerical example.
    The results confirm that topology can be analyzed using persistent homology and that the proposed selection operation significantly enhances the search performance of DDTD.
}




\maketitle

\section{Introduction}\label{sec1}

Structural optimization is a methodology aimed at maximizing desired performance by finding reasonable solutions through mathematical programming under computational models of physical phenomena.
Among these methodologies, topology optimization, initially proposed by \cite{bendsoe1988}, ensures maximum possible design freedom by designating material distribution within a given design domain as design variables.
Its potential to yield high-performance structures has led to a wide variety of engineering applications.

Currently, various topology optimization approaches have been developed, as summarized in \cite{sigmund2013} and \cite{deaton2014}.
Prominent examples include the density-based method \citep{bendsoe1989} and the level-set method \citep{allaire2004}, both of which update design variables based on sensitivity analysis, thus assuming differentiability of evaluation functions to be formulated as an optimization problem.
Moreover, in optimization problems with evaluation functions exhibiting strong multimodality, even if differentiable, extensive parameter studies are necessary to avoid convergence to undesirable local optima, but still do not always yield high-performance structures.
These challenges stemming from gradient-based optimizers pose significant barriers for further engineering applications, for example, minimax problems such as maximum stress minimization in stress-based topology optimization, and strongly multimodal optimization problems due to complex physics such as turbulence.

Focusing on optimization problems with non-differentiable or strongly multimodal evaluation functions, topology optimization with evolutionary algorithms (EAs) \citep{mitchell1999} has been developed.
EAs, as typified by genetic algorithms (GAs) \citep{goldberg}, are optimizers based on multi-point searching that mimic the emergent process of living organisms.
Various EA-based topology optimization methods have been proposed depending on the choice of the algorithm and the representation of design variables, which corresponds to the genotype in GAs.
\cite{chapman1994, wang2005}, \cite{madeira2010} and \cite{nimura2024} proposed GA-based methods with bit-array, graph, and quadtree representations, respectively.
\cite{wu2010} and \cite{luh2011} proposed other methods using differential evolution and particle swarm optimization with bit-array representation, respectively.
\cite{fujii2018} proposed another method using the covariance matrix adaption evolution strategy (CMA-ES) with level-set boundary representation.
Although these methods can yield reasonable material distributions as optimized solutions even for complex problems, EA-based topology optimization methods are typically challenged by the \textit{curse of dimensionality}.
This issue arises because the computational cost increases exponentially with the length of the design variables, or genetic sequences, limiting the dimensionality of optimization problems with an increasing number of design variables, i.e., the length of gene strings.
\cite{sigmund2011} has pointed out that the insufficient number of elements causes inaccurate physics and the loss of design freedom, resulting in only coarse optimized structures.

Data-driven design through the incorporation of machine learning offers a promising approach to avoiding the curse of dimensionality in EA-based frameworks.
As reviewed by \citep{woldseth2022}, data-driven topology optimization can be categorized into several approaches, among which the use of machine learning techniques for direct design, dimensionality reduction, and generative design is considered particularly effective.
As examples of direct design methods, \cite{yu2019} and \cite{behzadi2022} have proposed approaches using generative adversarial networks (GANs), a type of deep generative model.
These methods aim to predict optimal structures without iterations by leveraging pre-trained machine learning models.
As an example of dimensionality reduction, \cite{guo2018} have proposed a structural design method using a variational autoencoder (VAE).
This approach effectively solves topology optimization problems by exploiting the latent space of the VAE through a GA.
As an example of generative design, \cite{oh2019} have proposed a framework using a GAN.
Their method successfully obtains a diverse range of designs by iteratively learning data generated through topology optimization with GANs, starting from a limited set of existing designs.
Common challenges associated with such data-driven design methods include the need for a vast dataset to train the machine learning models in advance and the fact that the resulting designs often perform worse compared to those obtained through traditional topology optimization, as pointed out by \cite{woldseth2022}.

To achieve gradient-free optimization with a high degree of design freedom by integrating the EA and data-driven design, \cite{yamasaki2021} proposed a framework of data-driven topology design (DDTD) using a deep generative model.
Its core idea is that design candidates are iteratively updated by repeatedly selecting the superior ones from the dataset and generating new data by a deep generative model trained with them.
\cite{yaji2022} introduced an operation equivalent to mutation and systematically deriving promising initial data, employing the concept of multifidelity design \citep{yaji2020}.
\cite{kii2024} introduced a sampling method named latent crossover for deep generative models, positioning DDTD as a GA-based topology optimization framework.
Compared to typical EA-based methods, DDTD employs significantly lower dimensional genotype---latent variables encoded by a deep generative model---for the high-dimensional phenotype---i.e., discrete representations of material distributions, thereby avoiding the curse of dimensionality.

As mentioned above, while most studies on EA-based topology optimization focus on how to represent material distributions to implement efficient crossover and mutation, there are still challenges in selection operations.
Population diversity is crucial to prevent premature convergence and to achieve global search in evolutionary algorithms, and selection plays an important role in maintaining diversity \citep{goldberg}.
While typical approaches involve retaining inferior solutions in the population, such simple methods become challenging for maintaining diversity in multi-objective optimization problems \citep{li2015}.
In addition, \cite{tanabe2020, li2023} have pointed out that for strong nonlinear optimization problems, the complexity of the relationship between the design variable space and the objective function space makes it further difficult to maintain the population diversity.

\begin{figure*}[t]
    \centering
    \includegraphics[width=\textwidth]{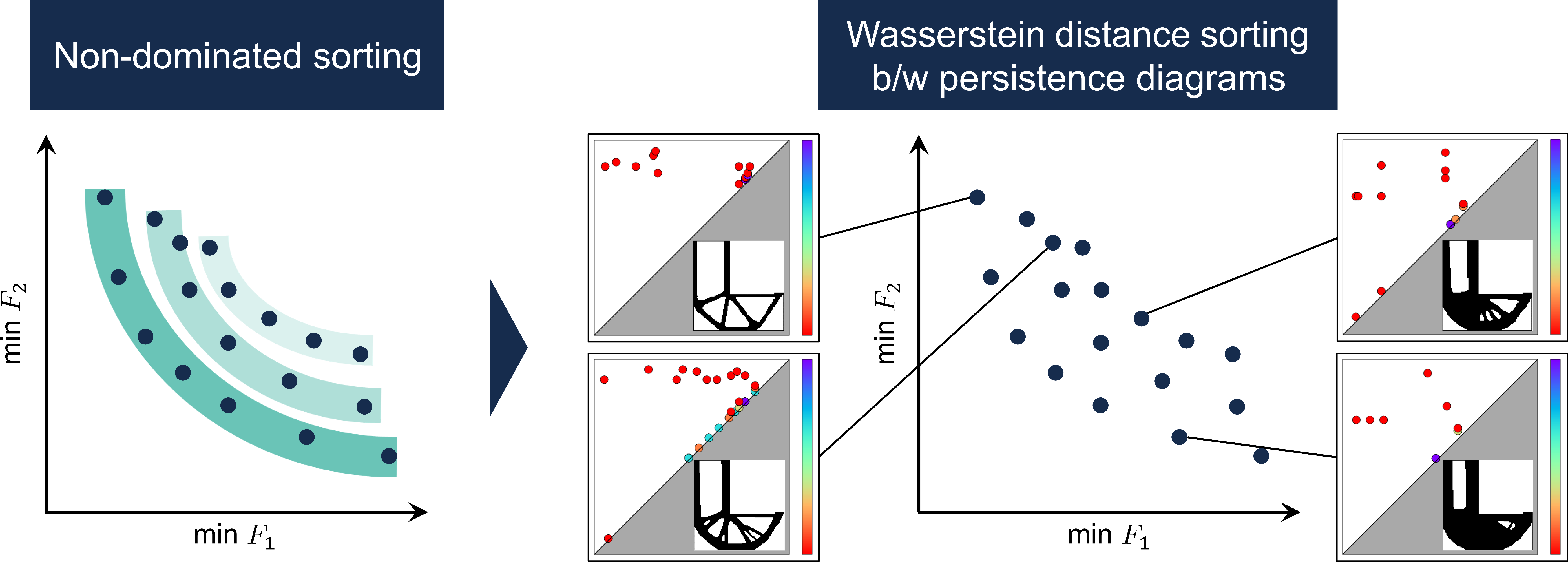}
    \caption{Schematic illustration of the proposed selection strategy}\label{fig_selection}
\end{figure*}

Given the above background, in this paper, we propose DDTD incorporating the selection operation to enhance the population diversity for multi-objective topology optimization problems with strong nonlinearity.
The key feature of the proposed selection strategy is its focus on selecting diverse solutions in the design variable space rather than the objective space.
However, quantifying the population diversity among material distributions represented by high-dimensional design variables is not straightforward. 
In this context, we consider that the differences between material distributions are characterized by the topological differences of the structures.
To capture these differences, we employ a topological data analysis method called persistent homology (PH) \citep{edelsbrunner2002, zomorodian2005} and incorporate a sorting based on the analyzed topological features into the selection process, as shown in Fig.~\ref{fig_selection}.
As a numerical example, we apply the proposed DDTD to a two-dimensional structural design problem and demonstrate the usefulness of PH for evaluating the topological features of material distribution data.
Through comparison with optimization results by the original DDTD, we verify the effectiveness of the proposed selection operation.

The rest of this paper is organized as follows.
In Section~\ref{sec2}, we describe the details of DDTD and discuss issues related to the selection operation as the motivation for using PH.
In Section~\ref{sec3}, we explain the proposed selection strategy in detail.
In Section~\ref{sec4}, the proposed method is applied to the structural design problem of an L-bracket, which is known as a benchmark for stress-based topology optimization, and the results are discussed.
Finally, Section~\ref{sec5} concludes the paper.

\section{Overview of data-driven topology design}\label{sec2}

\subsection{Optimization problem formulation}
Data-driven topology design (DDTD) \citep{yamasaki2021, yaji2022} is targeted at solving a multi-objective topology optimization problem formulated in the general form as follows:
\begin{equation}\label{eq_moto}
    \begin{aligned}
        & \underset{\rho}{\text{minimize}}
            && \left[F_{1}(\rho), F_{2}(\rho), \ldots , F_{r_{\text{o}}}(\rho)\right],\\
        & \text{subject to}
            && G_{j}(\rho)\leq0 \quad \text{for } j=1, 2, \ldots, r_{\text{c}},\\
            &&& \rho (\boldsymbol{x})\in\{0, 1\},\quad\forall \boldsymbol{x}\in D,
    \end{aligned}
\end{equation}
where $F_{i}\,(i=1, 2, \ldots, r_{\text{o}})$ and $G_{j}\,(j=1, 2, \ldots, r_{\text{c}})$ are the objective and constraint functions, respectively.
Material distributions are represented as the design variable $\rho (\boldsymbol{x})$, where $\boldsymbol{x}$ is the coordinates at an arbitrary point within the fixed design domain $D$.
The design variable $\rho$ takes discrete values of 0 or 1, where $\rho (\boldsymbol{x})=1$ and $0$ mean the material and void, respectively.
The original topology optimization problem in \eqref{eq_moto} is often difficult to solve directly because it involves a highly nonlinear optimization problem with significant multimodality and can be formulated using non-differentiable evaluation functions.
Therefore, for some procedures described later, we formulate the low-fidelity optimization problem using the idea of multifidelity design \citep{yaji2020} as follows:
\begin{equation}\label{eq_lf}
	\begin{aligned}
		& \underset{\rho^{(k)}}{\text{minimize}}
			&& \widetilde{F}(\rho^{(k)}),\\
		& \text{subject to} 
			&& \widetilde{G}_{l}(\rho^{(k)}, \boldsymbol{s}^{(k)})\leq0 \quad \text{for } l=1, 2, \ldots, \widetilde{r}_{\text{c}},\\
			&&& \rho^{(k)}(\boldsymbol{x})\in[0, 1],\quad\forall \boldsymbol{x}\in D,\\
		& \text{for given}
			&& \boldsymbol{s}^{(k)},
	\end{aligned}
\end{equation}
where $\widetilde{F}$ and $\widetilde{G}_{l}\,(l=1, 2, \ldots, \widetilde{r}_{\text{c}})$ are the objective and constraint functions of a low-fidelity optimization problem, respectively.
They are simplified pseudo-functions for the original ones $F_{i}$ and $G_{j}$, formulated to be computationally easier and differentiable.
The low-fidelity optimization problem in \eqref{eq_lf} is assumed to be solved by typical topology optimization methods such as the density-based method \citep{bendsoe2003}, and is reformulated as a single-objective optimization problem on the basis of the $\varepsilon$-constraint method \citep{haimes1971} or the weighted-sum method \citep{zadeh1963}, as opposed to the original multi-objective problem in \eqref{eq_moto}.
Additionally, the design variable $\rho^{(k)}$ is relaxed to continuous values in the range of 0 to 1. 
$\boldsymbol{s}=[s_{1}, s_{2}, \ldots, s_{N_{\text{sd}}}]$ represents the set of $N_{\text{sd}}$ types of artificial design parameters called seeding parameters, and $\boldsymbol{s}^{(k)}$ represents the sample point of $\boldsymbol{s}$.
The seeding parameter $\boldsymbol{s}$ includes optimization parameters such as filter radius and projection method parameters in density-based optimization, as well as constraint values.

\subsection{Optimization procedure}
\begin{figure*}[t]
    \centering
    \includegraphics[width=0.9\textwidth]{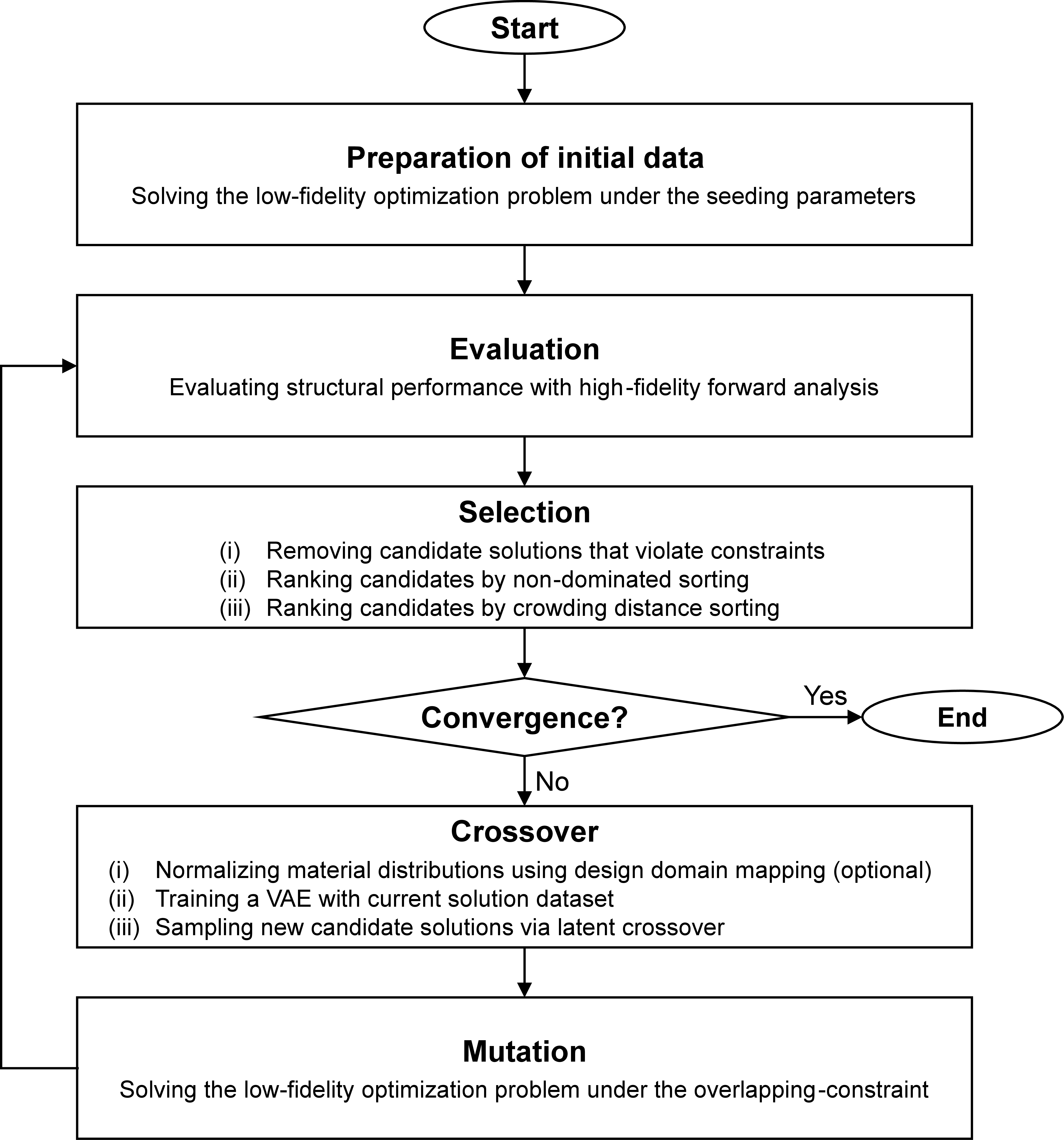}
    \caption{Optimization flowchart of DDTD}\label{fig_flowchart}
\end{figure*}

After formulating the optimization problems as described above in advance, DDTD performs topology optimization through genetic algorithm-based procedures, i.e., iteratively updating solutions through three genetic operations: selection, crossover, and mutation.
Figure~\ref{fig_flowchart} outlines the optimization flowchart of DDTD, and the details of each procedure are explained here.

\subsubsection{Preparation of initial data set}
Diverse and promising initial solutions for the original problem in \eqref{eq_moto} are derived by solving a simplified problem, the low-fidelity optimization problem in \eqref{eq_lf}, under various parameter settings.

\subsubsection{Evaluation}
The candidate solutions are evaluated by the high-fidelity analysis model with the objective function $F_{i}$ and constraint functions $G_{j}$ in the original problem in \eqref{eq_moto}.
Note that design variables are binarized to $\{0, 1\}$ for high-fidelity evaluation, and only forward analysis for the evaluation functions $F_{i}$ and $G_{j}$ is required.

\subsubsection{Selection}
Based on the objective values from high-fidelity evaluation, superior candidate solutions are selected to be preserved for the next generation.
Since the target problem in \eqref{eq_moto} is a multi-objective optimization problem, it is necessary to construct the dataset of solution sets under Pareto optimality.
Here, the selection process in DDTD is responsible for constructing the VAE training data.
Note the distinction from a typical GA, where selection involves not only eliminating inferior individuals but also choosing parents for crossover.
The details of the selection operation in DDTD are discussed in Section~\ref{sec23}.

\subsubsection{Convergence check}
The optimization computation is checked for convergence.
When either the pre-determined maximum number of iterations is reached or the hypervolume indicator \citep{shang2020}, a convergence performance measure in multi-objective optimization, has converged sufficiently, the optimization computation is terminated.

\subsubsection{Crossover}
\begin{figure*}[t]
    \centering
    \includegraphics[width=0.8\textwidth]{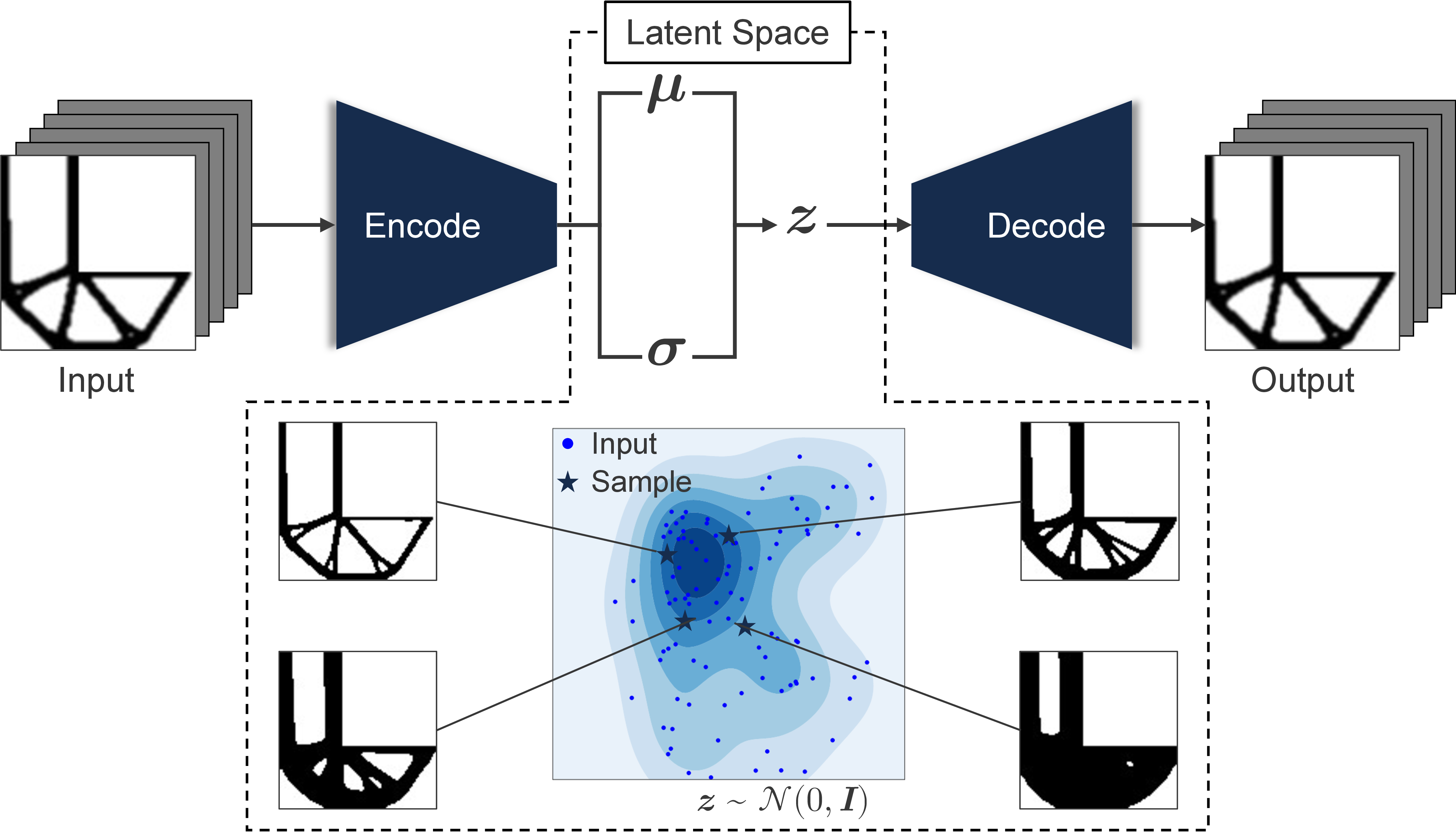}
    \caption{Schematic illustration of VAE}\label{fig_vae}
\end{figure*}

A deep generative model is trained with the dataset of solution sets constructed in the selection process.
Representative deep generative models include variational autoencoders (VAEs) \citep{kingma2013} and generative adversarial networks (GANs) \citep{goodfellow2014}, and in prior studies of DDTD \citep{yamasaki2021, yaji2022, kii2024}, VAEs have been preferred due to their learning stability.
As shown in Fig.~\ref{fig_vae}, a VAE consists of two neural networks called the encoder and the decoder; the former compresses high-dimensional input data to low-dimensional latent variables $\boldsymbol{z}$ with their mean $\boldsymbol{\mu}$ and standard deviation $\boldsymbol{\sigma}$, while the latter reconstructs the original dimensional output data from latent variables.
A Gaussian distribution is assumed in the latent space of a VAE by defining the latent variable $\boldsymbol{z}$ as follows:
\begin{equation}\label{eq_z}
    \boldsymbol{z}=\boldsymbol{\mu}+\boldsymbol{\sigma}\odot\boldsymbol{\varepsilon},
\end{equation}
where $\odot$ represents the element-wise product and $\boldsymbol{\varepsilon}$ is a random sample from a Gaussian distribution $\mathcal{N}(0, \boldsymbol{I})$.
Based on the definition of latent variables in \eqref{eq_z}, the VAE assumes a probability distribution in the latent space and functions as a generative model, generating new data that inherits the features of the training data through sampling from the latent space and reconstructing it using the decoder.
To define latent variables as the genotype of EAs in DDTD, latent crossover proposed by \cite{kii2024} is employed as a sampling technique.
Note that in the original DDTD papers \citep{yamasaki2021, yaji2022}, design domain mapping \citep{yamasaki2019} is employed to map material distributions into the unit square domain for normalization of training data.
This normalization step is not always necessary for VAE training and is provided as an option.

\subsubsection{Mutation}
The low-fidelity optimization problem in \eqref{eq_lf} is solved under the following overlapping-constraint:
\begin{equation}\label{eq_mut}
	\widetilde{G}_{\text{mut}}(\rho^{(m)})=\int_{D}\rho^{(m)}\rho^{\text{ref}}\,\text{d}\Omega- G_{\text{mut}}^{\text{max}}\int_{D}\,\text{d}\Omega\leq 0,
\end{equation}
where $m=1, 2, \ldots, N_{\text{mut}}$ is the index of mutants and $G_{\text{mut}}^{\text{max}}$ is a parameter that controls the degree of overlapping between design variables $\rho^{(m)}$ and the reference density distribution $\rho^\text{ref}$, represented as follows:
\begin{equation}\label{eq_ref}
	\rho^{\text{ref}}(\boldsymbol{x})=\frac{1}{N_{\text{pop}}} \sum_{n=1}^{N_{\text{pop}}}\rho^{(n)}(\boldsymbol{x}).
\end{equation}
Herein, $N_{\text{pop}}$ is the population size, and the reference density distribution in \eqref{eq_ref} denotes the superposition of material distributions of all solutions in that generation, i.e., the average material distribution.
In this way, by solving an easily solvable pseudo-problem, the low-fidelity optimization problem in \eqref{eq_lf}, under the overlapping-constraint (3) with the average material distribution (4), promising material distributions with new features not present in the current solution set are injected as mutants.

\subsection{Issues in selection}\label{sec23}
DDTD solves multi-objective topology optimization problems and uses the selection strategy of the non-dominated sorting genetic algorithm II (NSGA-II) \citep{deb2002, verma2021}, one of the representative GAs for multi-objective problems.
The details of its selection operation are provided in Section~\ref{sec231}, and challenges in solving topology optimization problems are discussed in Section~\ref{sec232}.

\subsubsection{Selection operation of NSGA-II}\label{sec231}
\begin{figure}[t]
    \centering
    \includegraphics[width=0.95\textwidth]{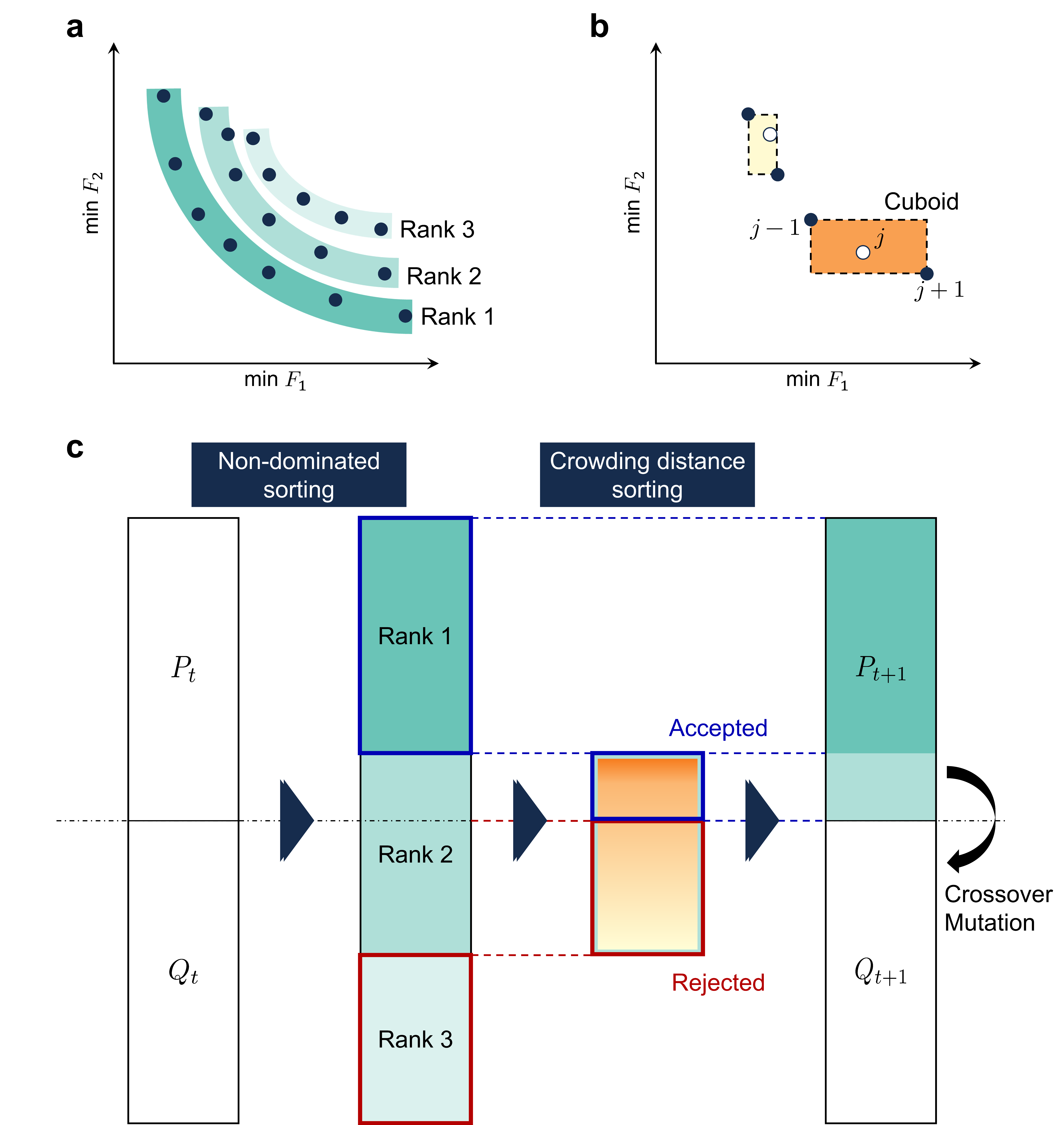}
    \caption{Schematic illustration of selection operation in NSGA-II: (a) non-dominated sorting; (b) crowding distance sorting; (c) overall procedure}\label{fig_nsga}
\end{figure}

The NSGA-II selection operation consists of ranking by two sortings: non-dominated sorting and crowding distance sorting.
The characteristics of the rules for selecting one out of the two candidates in the NSGA-II procedure can be summarized in the following two points \citep{verma2021}:
\begin{enumerate}
    \item If two candidates have different ranks in the non-dominated sorting, then the one with the better rank is selected for the next generation.
    \item If two candidates have the same ranks in the non-dominated sorting, the one with the larger crowding distance is selected for the next generation.
\end{enumerate}
An outline of such selection rules is illustrated in Fig.~\ref{fig_nsga}, and their details are described here.

\textbf{Non-dominated sorting}
In this procedure, candidate solutions are sorted based on the concept of Pareto dominance.
Here, for the multi-objective minimization problem in \eqref{eq_moto}, the candidate $\rho ^{(2)}$ is defined as dominating the candidate $\rho ^{(1)}$ when the following conditions hold:
\begin{equation}\label{eq_nd}
    F_i(\rho ^{(1)})\leq F_i(\rho ^{(2)})~(\forall i=1, 2, \ldots, r_{\text{o}}),\quad F_i(\rho ^{(1)})< F_i(\rho ^{(2)}) ~(\exists i=1, 2, \ldots, r_{\text{o}}).
\end{equation}
Based on the definition in \eqref{eq_nd}, non-dominated solutions that are not dominated by any other candidates are initially assigned rank 1 and extracted from the population to form the first front. 
Subsequently, rank 2 candidates are assigned similarly.
As shown in Fig.~\ref{fig_nsga}a, This iterative process continues until all candidate solutions are assigned to fronts.

\textbf{Crowding distance sorting}
In this procedure, the crowding degree of the front in the objective space is calculated for candidate solutions assigned the same rank during the non-dominated sorting process.
As shown in Fig.~\ref{fig_nsga}b, the crowding distance $d^{(j)}$ for the $j$-th candidate $\rho ^{(j)}$ is calculated based on the cuboid formed by neighboring ones as follows:
\begin{equation}\label{eq_cd}
    d^{(j)}=\sum_{i=1}^{r_{\text{o}}}\frac{F_i(\rho ^{(j+1)})-F_i(\rho ^{(j-1)})}{F_i^{\text{max}}-F_i^{\text{min}}},
\end{equation}
where $F_i^{\text{max}}$ and $F_i^{\text{min}}$ are the maximum and minimum values of objective function $F_i$ among all the candidates, respectively.

An overview of the overall NSGA-II procedure with the above two sorts is shown in Fig.~\ref{fig_nsga}c.
Here we consider selecting a population $P_{t+1}$ of size $N_{\text{pop}}$ from the dataset $P_t\cup Q_t$, in which $P_t$ and $Q_t$ are the current population consisting of the solution set at generation $t$ and new candidate solutions obtained through crossover and mutation, respectively.
If the number of rank 1 candidates is greater than or equal to $N_{\text{pop}}$, then $N_{\text{pop}}$ solutions with larger 
crowding distances in \eqref{eq_cd} are selected from the first front to form $P_{t+1}$.
In contrast, if the number of rank 1 candidates is less than $N_{\text{pop}}$, the first front is directly transferred to $P_{t+1}$, and the remaining solutions are selected from the second front with larger crowding distances and added to $P_{t+1}$.
If the size of $P_{t+1}$ is still less than $N_{\text{pop}}$, this process is repeated with the subsequent fronts until the size of $P_{t+1}$ reaches $N_{\text{pop}}$.
The constructed $P_{t+1}$ undergoes crossover and mutation to form a new set of candidate solutions $Q_{t+1}$ for the next generation $t+1$, and the same selection process is applied to the new dataset $P_{t+1}\cup Q_{t+1}$ to form the next population $P_{t+2}$.
In the optimization process of NSGA-II, this series of genetic operations is repeated until the convergence criteria are satisfied.

\subsubsection{Challenges in maintaining diversity}\label{sec232}
Selection plays an important role in EAs for maintaining population diversity and facilitating global search, and various selection schemes have been proposed, such as roulette wheel selection, ranking selection, and tournament selection in GAs \citep{goldberg}.
These approaches involve selecting not only superior solutions with higher fitness values but also inferior solutions to preserve a diverse set of solutions for the next generation.
This strategy is also exemplified by NSGA-II discussed in Section~\ref{sec231}, which employs the crowding distance sorting based on a similar principle.
In other words, DDTD considers the diversity of material distribution $\rho$ based on the objectives $F_i$ in the multi-objective topology optimization problem in \eqref{eq_moto}.

\begin{figure}[t]
    \centering
    \includegraphics[width=0.8\textwidth]{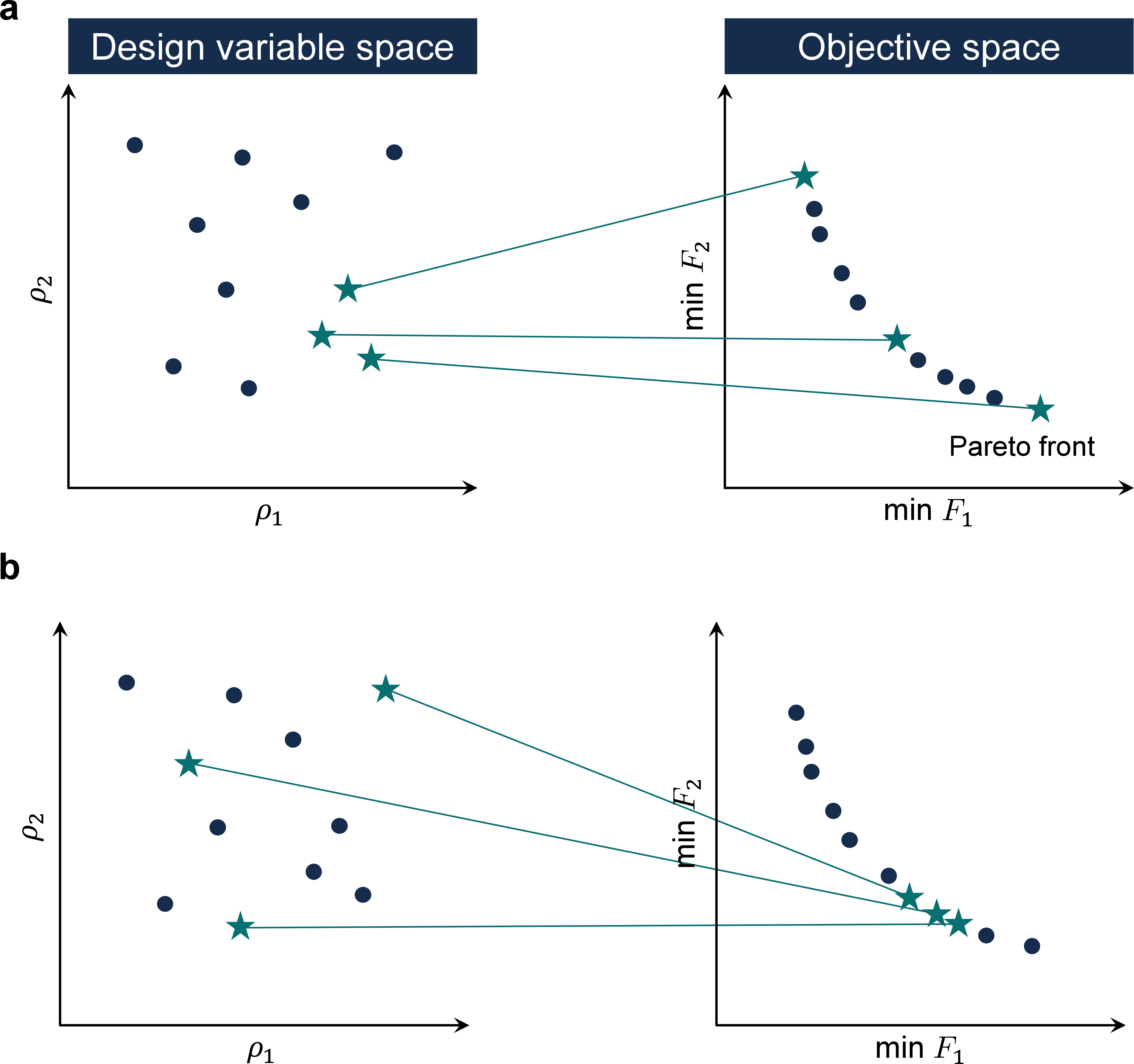}
    \caption{Relationship between design variable space and objective space in multi-objective optimization problems with significant nonlinearity: (a) Solutions with similar design variables for different objective values; (b) Solutions with different design variables for similar objective values}\label{fig_ds_os}
\end{figure}

\cite{tanabe2020, li2023} have pointed out an issue with such selection approaches in multi-objective optimization problems with strong nonlinearity.
In the case of strong nonlinear optimization problems, the evaluation functions are extremely sensitive to the design variables, i.e., even slight changes in the variables may cause significant fluctuations in the objective values.
As a result, even if distant solutions are selected in the objective space, they may still be located in proximity in the design variable space, as shown in Fig.~\ref{fig_ds_os}a, which shows an example of a two-objective optimization problem with two design variables for simplicity.
The reverse is also true, as shown in Fig.~\ref{fig_ds_os}b, where solutions that are close in the objective space may be far apart in the design variable space.
In other words, conventional selection approaches, which ensure diversity based on objective values, may result in a population filled with solutions with similar design variables, and crossover and mutation cannot produce new candidates, potentially leading to premature convergence.
Particularly in topology optimization, where complex physics often leads to strongly nonlinear optimization problems, it becomes crucial to adopt a selection strategy that maintains the intrinsic population diversity in the design variable space for global search through EAs.

\begin{figure}[t]
    \centering
    \includegraphics[width=0.8\textwidth]{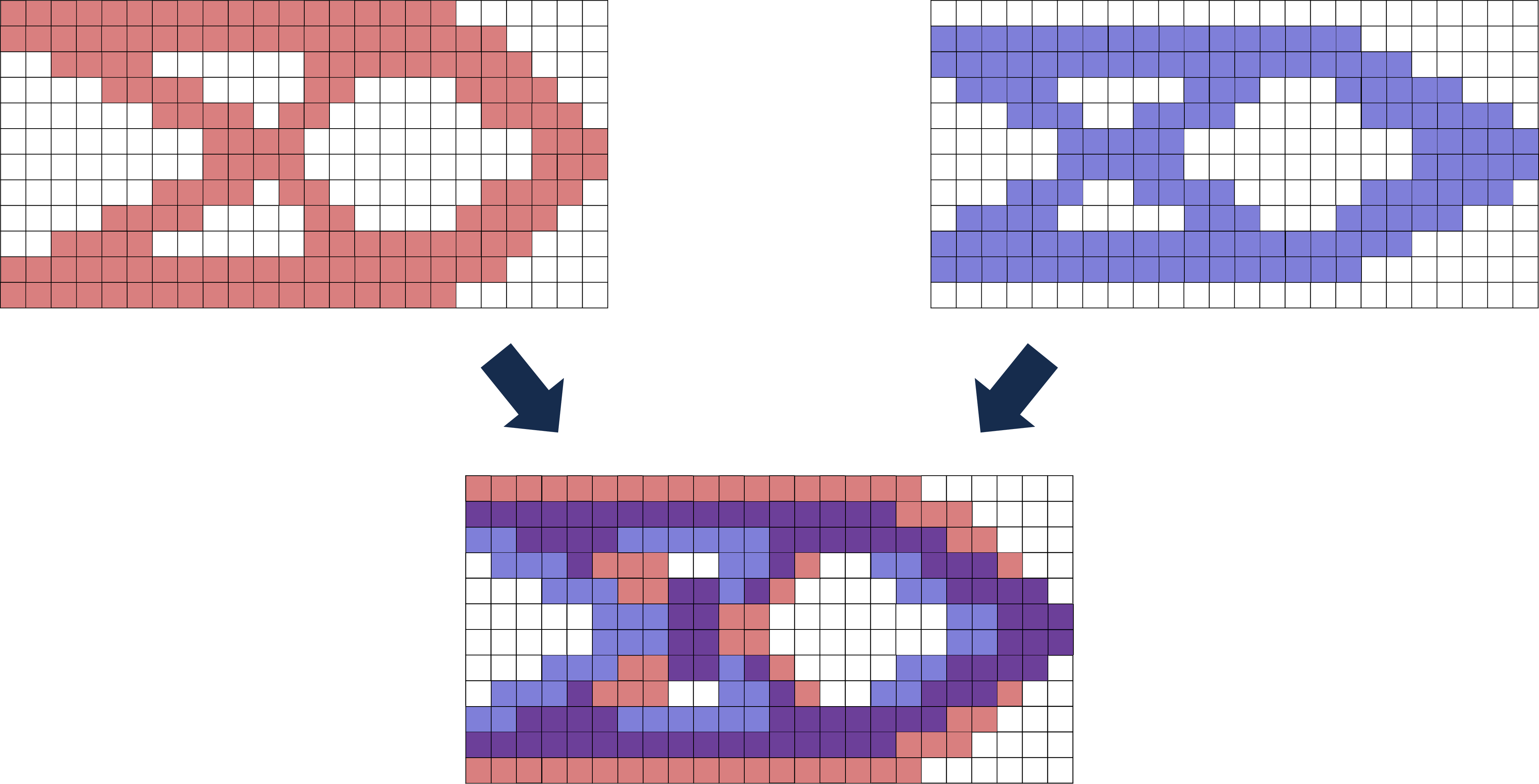}
    \caption{Examples of material distributions where structural differences are difficult to quantify in the $L^p$ norm}\label{fig_lp}
\end{figure}

It is theoretically possible to quantitatively measure diversity in the design variable space, i.e., to quantify the differences between material distribution in topology optimization and incorporate this into the selection algorithm.
For the topology optimization problem in \eqref{eq_moto} to be solved, the simplest method to measure the difference between material distributions $\rho ^{(1)}$ and $\rho ^{(2)}$ is the $L^p$ norm, defined as follows:
\begin{equation}\label{eq_lp}
	\left\lVert \rho ^{(1)}-\rho ^{(2)}\right\rVert _p= \left(\int_{D}\left\lvert \rho ^{(1)}-\rho ^{(2)}\right\rvert ^p\,\text{d}\Omega\right) ^{\frac{1}{p}}.
\end{equation}
While such norm measures are commonly used to evaluate the distance between functions or vectors, their effectiveness in capturing the structural differences between material distributions in topology optimization is questionable.
Although the material distribution $\rho (\boldsymbol{x})$ in the optimization problem in \eqref{eq_moto} is theoretically represented as a continuous function, in typical topology optimization such as the density-based method \citep{bendsoe2003}, it is generally discretized into a finite number of design variables for computational procedures.
The number of elements required for discretization typically exceeds several thousand to achieve a high degree of design freedom.
Even if the norm in \eqref{eq_lp} is calculated for such a design variable vectors, it merely sums up the differences of each element and does not measure the essential structural differences between the material distributions.
For example, the $L^p$ norm between structures where each member is shifted by one pixel, as shown in Fig.~\ref{fig_lp}, would be excessively large.
A typical example of a distance function that resolves these $L^p$ norm problems is the Wasserstein metric based on the idea of optimal transport.
However, computing the Wasserstein distance requires solving a minimization problem about the transportation cost of material from one distribution to another within the design domain.
This process becomes computationally impractical for thousands of discretized material distributions, as it would require repeatedly solving this optimization problem during each selection operation.

\section{Proposed selection strategy}\label{sec3}

To overcome the challenges of selection in DDTD discussed in Section~\ref{sec232}, it is necessary to employ an efficient diversity measure of material distributions with low computational cost.
In structural optimization problems, the overall shape of a structure is typically determined uniquely by the predefined design domain and boundary conditions, suggesting that the diversity among structures is largely influenced by their topology.
Therefore, this paper focuses on a topological data analysis method as a means to condense the topological information of material distribution data, based on the premise that topology is the crucial determinant of population diversity.
Then, by applying a measure based on optimal transport to calculate the distance on the condensed topological data, the differences between material distributions can be quantified at a practical computational cost. 
Specifically, the topological features extracted using a topological data analysis method called persistent homology \citep{edelsbrunner2002, zomorodian2005} are represented as persistence diagrams, and the distance between them is quantified by the Wasserstein distance \citep{mileyko2011} and incorporated into the selection operation.
In the following, an overview of persistent homology and the proposed selection operation are described in detail.

\subsection{Persistent homology}
\subsubsection{Overview}\label{sec311}
\begin{figure}[t]
    \centering
    \includegraphics[width=0.95\textwidth]{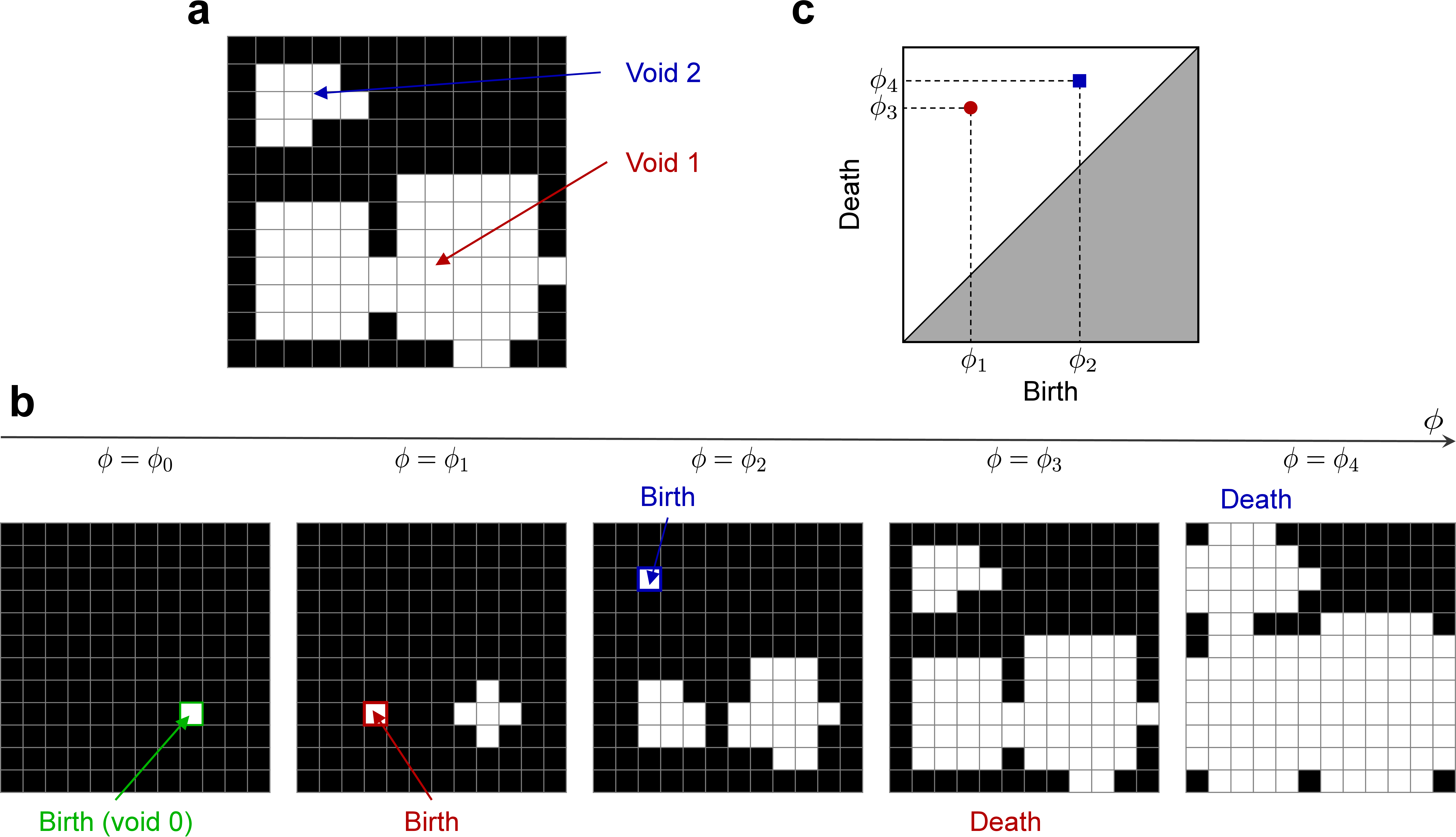}
    \caption{Schematic illustration of persistent homology: (a) target binary image; (b) filtration based on signed distance function with the Manhattan distance; (c) persistence diagram}\label{fig_ph}
\end{figure}

Persistent homology (PH) \citep{edelsbrunner2002, zomorodian2005} is a type of topological data analysis method that mathematically quantifies geometric features of targeting complex data using the concept of topology.
Its scope of application covers a wide variety of data sets, including point clouds, images, graphs, and so on.
Here, topology refers to geometric properties and spatial relations unaffected by the continuous change of shape or size, which is different from its context in structural optimization.
The fundamental idea of PH is to track the birth and death of topological features over time through an operation called filtration, which involves gradually increasing a scale parameter \citep{munch2017, otter2017}.
Assuming an application to material distribution data in topology optimization, PH for a binary image is shown in Fig.~\ref{fig_ph}.
Here, we discuss 0th persistent homology, which involves identifying voids by considering the connected components of the white regions in Fig.~\ref{fig_ph}a.
First, we consider a level-set function $\phi$ that assigns integer values to each pixel using a signed distance function with the Manhattan distance for the boundary between white and black pixels \citep{obayashi2018}.
The values of $\phi$ serve as the scale parameter in this example, and the filtration considers shapes obtained as their union set.
In Fig.~\ref{fig_ph}a, two voids can be visually identified in the binary image under analysis, and in the filtration process in Fig.~\ref{fig_ph}b, as $\phi$ increases from its minimum value, they appear at a certain stage and eventually disappear by merging with other voids in the resulting shape.
More specifically, void 0 first appears at $\phi=\phi_0$, and void 1, which appears at $\phi=\phi_1$, and disappears by merging with void 0 at $\phi=\phi_2$.
Similarly, void 2 is born at $\phi=\phi_3$ and dies at $\phi=\phi_4$.
These birth-death pairs are represented as points on a two-dimensional plot, known as a persistence diagram (PD), with the birth time on one axis and the death time on the other.
The persistence diagram $D$ can be defined as follows:
\begin{equation}
    D = \{(b_i, d_i): i=1,2,\dots, n\},
\end{equation}
where $(b_i, d_i)$ is the $i$-th birth-death pair and $n$ is the number of birth-death pairs, respectively.
For example, the PD in Fig.~\ref{fig_ph}c is given by $\{(\phi_1, \phi_3), (\phi_2, \phi_4)\}$, indicating that the two voids corresponding to the two plots on the PD are present in the original binary image in Fig.~\ref{fig_ph}a.
Note that the birth-death pair for void 0, which appears first and persists as $\phi$ increases, is $(\phi_0, \infty)$, and it is not included in the PD.
The value $d_i-b_i$, indicating the length of time a feature persists from its appearance to its disappearance, is called the lifetime, which corresponds to how far the plot is from the diagonal.
In the persistence diagram of Fig.~\ref{fig_ph}c, the plot corresponding to void 1 is further from the diagonal than that of void 2, indicating that void 1 has a greater persistence with a larger lifetime.
In this way, the sizes, numbers, and spatial relationships of connected components, holes, and voids analyzed by PH are summarized in the PD, enabling the rational extraction of these topological features from the targeting data set.

\subsubsection{Quantifying topological difference}\label{sec312}
\begin{figure}[t]
    \centering
    \includegraphics[width=\textwidth]{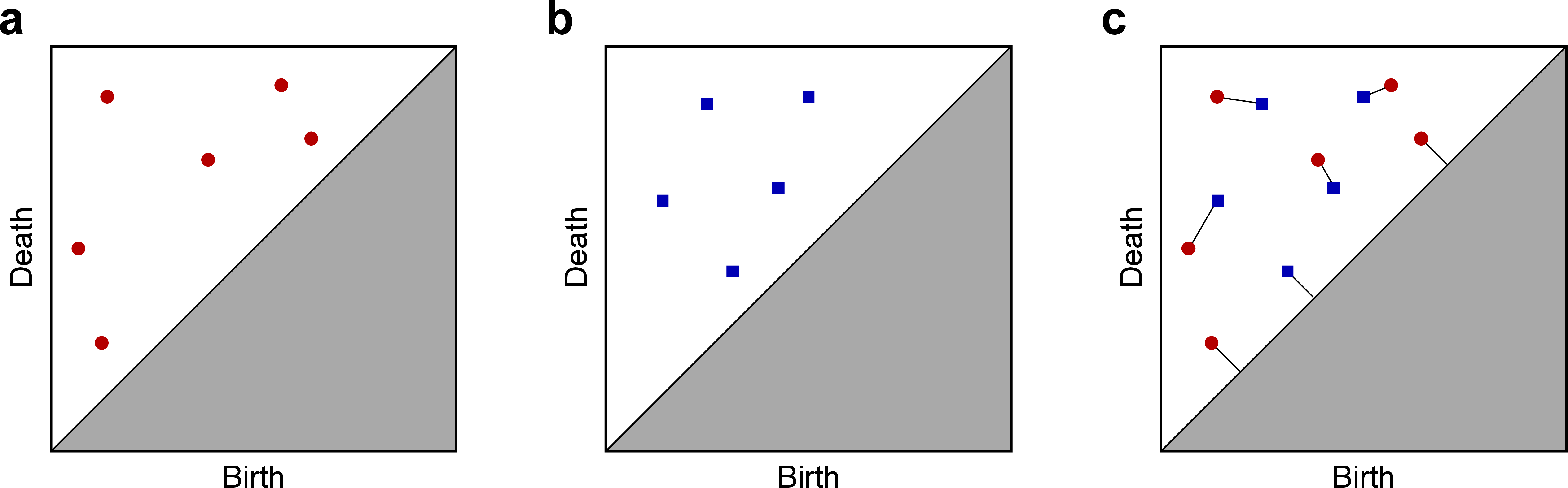}
    \caption{Schematic illustration of partial matching between persistence diagrams (PDs) that minimizes the transportation cost in Wasserstein distance calculation: (a) PD1; (b) PD2; (c) partial matching between PD1 and PD2}\label{fig_pds}
\end{figure}

To compare complex data using PH, distance metrics between PDs have been proposed \citep{mileyko2011}.
The basic idea is based on the concept of optimal transport, where the matching of points between two PDs is considered.
The cost of the matching, minimized to the least possible value, is used as the distance between the two PDs.
However, based on the concept of PH discussed in Section~\ref{sec311}, points near the diagonal of PDs correspond to holes that disappear as soon as they appear, and they are insignificant points like noisy plots and should not affect the matching cost.
Therefore, let $q_1\in D_1$ and $q_2\in D_2$ be the points on $D_1$ and $D_2$, and the partial matching between $D_1$ and $D_2$ is given by a subset $M\subset D_1\times D_2$ as follows:
\begin{itemize}
    \item For $\forall q_1\in D_1$, there is at most one $q_2\in D_2$ such that $(q_1, q_2)\in M$.
    \item For $\forall q_2\in D_2$, there is at most one $q_1\in D_1$ such that $(q_1, q_2)\in M$.
\end{itemize}
Such a partial matching $M$ is represented as $M: D_1\leftrightarrow D_2$, where $(q_1, q_2)\in M$ denotes a pair of matched points on $D_1$ and $D_2$.
On the other hand, the remaining unmatched points are denoted as $q\in D_1\sqcup  D_2$, then the transportation cost based on the concept of optimal transport is formulated as follows:
\begin{equation}\label{eq_cost}
    c_p(M)=\left(\sum_{(q_1, q_2)\in M} (\left\lVert q_1 - q_2 \right\rVert _p)^p + \sum_{q \in D_1 \sqcup D_2 } (\left\lVert q - \pi(q) \right\rVert _p)^p \right)^{\frac{1}{p}},
\end{equation}
where $\pi(q)$ is the orthogonal projection of $q$ onto the diagonal.
The Wasserstein distance between $D_1$ and $D_2$ is calculated as the minimum transportation cost of \eqref{eq_cost} as follows:
\begin{equation}\label{eq_w}
    W_p(D_1, D_2)=\underset{M\colon D_1\leftrightarrow D_2}{\inf} c_p(M).
\end{equation}
An example of partial matching between persistence diagrams in the calculation of Wasserstein distance is shown in Fig.~\ref{fig_pds}.
Figures~\ref{fig_pds}a and \ref{fig_pds}b are PDs with 6 and 5 points, respectively, it is not possible to consider a point-to-point matching for all plots due to the differing number of them.
The transportation cost in \eqref{eq_cost} allows us to consider matching not only to points but also to the diagonal, resulting in the optimal partial matching as shown in Fig.~\ref{fig_pds}c.
Additionally, in the calculation of the Wasserstein distance using the transportation cost in \eqref{eq_cost}, noise-like plots with nearly the same birth and death times are matched with the diagonal and have little effect on the cost.
The commonly used Wasserstein distance $W_2$ at $p=2$ in \eqref{eq_w} corresponds to the sum of the lengths of the line segments connecting the matchings in Fig.~\ref{fig_pds}c.

\subsubsection{Previous research and novelty in this study}\label{sec313}
Persistent homology, which enables the analysis of topological features, has a high affinity with topology optimization, which targets the shape and phase of a structure, and a few previous studies have been reported in recent years.
\cite{wang2022} proposed a method for implementing topological constraints in the SIMP method \citep{bendsoe1999} based on the concept of PH.
They successfully obtained a compliance minimization design that satisfies inequality constraints related to the holes within the structure.
\cite{depeng2024} proposed a method to determine the effective relative density range of triply periodic minimal surfaces (TPMSs) based on PH.
They successfully obtained high-stiffness porous structures through topology optimization by determining the effective thresholds of TPMSs from a topological perspective using PH.
\cite{behzadi2022, hu2024} proposed a topology optimization framework using generative models, specifically GANs and VAEs, respectively, incorporating a loss function based on PH.
They showed that training the neural network with the distance between persistence diagrams as a topological loss improves the connectivity of the generated structures compared to general generative models that minimize only the reconstruction loss.

In this paper, we propose a selection operation of EAs that enhances population diversity, focusing on the quantification of differences in material distributions as the Wasserstein distance between PDs.
In other words, compared to the aforementioned previous studies, the novelty of this research lies in the topology optimization framework with an EA that incorporates topological features analyzed using PH.
In this way, topological differences between complex data can be quantified as Wasserstein distances between PDs.

\subsection{Wasserstein distance sorting between persistence diagrams}\label{sec32}
Based on the selection operation of NSGA-II described in Section~\ref{sec231}, this paper proposes a selection strategy that incorporates the analysis of topological features of material distributions using PH. 
Specifically, to address the challenges of directly using the NSGA-II selection operation in DDTD described in Section~\ref{sec232}, we propose a new sorting method, named Wasserstein distance sorting between PDs, as an alternative to crowding distance sorting.
The details of the proposed sorting procedure are as follows:
\begin{enumerate}
    \item For all candidate solutions, PH is computed from the material distribution $\rho ^{(i)}$ to obtain the corresponding PD, $D^{(i)}$.
    \item For each $D^{(i)}$, the Wasserstein distance is calculated in a pairwise manner to generate the following $N_{\text{cand}}\times N_{\text{cand}}$ distance matrix:
        \begin{equation}
            A=
            \begin{pmatrix}
                W_p(D^{(1)}, D^{(1)})               & W_p(D^{(1)}, D^{(2)})                 & \cdots & W_p(D^{(1)}, D^{(N_{\text{cand}})}) \\
                W_p(D^{(2)}, D^{(1)})               & W_p(D^{(2)}, D^{(2)})                 & \cdots & W_p(D^{(2)}, D^{(N_{\text{cand}})}) \\
                \vdots                              & \vdots                                & \ddots & \vdots \\
                W_p(D^{(N_{\text{cand}})}, D^{(1)}) & W_p(D^{(N_{\text{cand}})}, D^{(2)})   & \cdots & W_p(D^{(N_{\text{cand}})}, D^{(N_{\text{cand}})})
            \end{pmatrix},
        \end{equation}
    where $N_{\text{cand}}$ is the total number of candidate solutions given by $N_{\text{cand}}=N_{\text{pop}}+N_{\text{VAE}}+N_{\text{mut}}$, where $N_{\text{cand}}$, $N_{\text{pop}}$ and $N_{\text{mut}}$ are the population size, the number of generated data by crossover using a VAE, and the number of mutants, respectively.
    \item The Wasserstein distance for the PD of the $i$-th candidate is calculated as the sum of the $i$-th row of the distance matrix $A$ as follows:
        \begin{equation}\label{eq_d}
            d^{(i)}=\sum_{j=1}^{N_{\text{cand}}}W_p(D^{(i)}, D^{(j)}).
        \end{equation}
    \item Based on the $d^{(i)}$ of \eqref{eq_d}, the candidate solutions are sorted in descending order.
\end{enumerate}

The sorting obtained through the above operations replaces the crowding distance sorting in NSGA-II.
The overall scheme of the proposed selection method is similar to that of NSGA-II: the non-dominated sorting ranks $N_{\text{cand}}$ candidate solutions based on Pareto dominance, while the Wasserstein distance sorting between PDs determines priority among the rank containing the $N_{\text{pop}}$-th candidate.
Note that if rank 1 is assigned to all $N_{\text{cand}}$ candidates in the non-dominated sorting, indicating that the optimization process has entered the convergence stage, the crowding distance sorting is employed instead of the proposed sorting to obtain a continuous Pareto front.
This division of the optimization process is based on exploration and exploitation in EAs \citep{crepinsek2013}, where the former corresponds to the Wasserstein distance sorting and the latter corresponds to the crowding distance sorting.

The proposed sorting method allows selection of GAs to preserve diverse solutions for the next generation based on the topology calculated from the material distribution using PH, whereas the conventional method selects them based on the objective values.
Structural topology is the most distinctive factor characterizing the diversity of material distributions in topology optimization.
It is expected that mating through crossover and mutation of GAs for topologically diverse material distribution allows the population to spread out more widely in the solution space, facilitating global search.

\section{Numerical example}\label{sec4}
In this section, we demonstrate the usefulness of the proposed selection strategy incorporated into DDTD through a numerical example.
First, we confirm that topology can be analyzed through the application of PH to material distributions and that the Wasserstein distances between persistence diagrams can be calculated appropriately.
Then, we verify the effectiveness of the proposed selection strategy by comparing the optimization results with those from the original DDTD.

\subsection{Problem settings}\label{sec41}
\begin{figure}[t]
    \centering
    \includegraphics[width=0.47\textwidth]{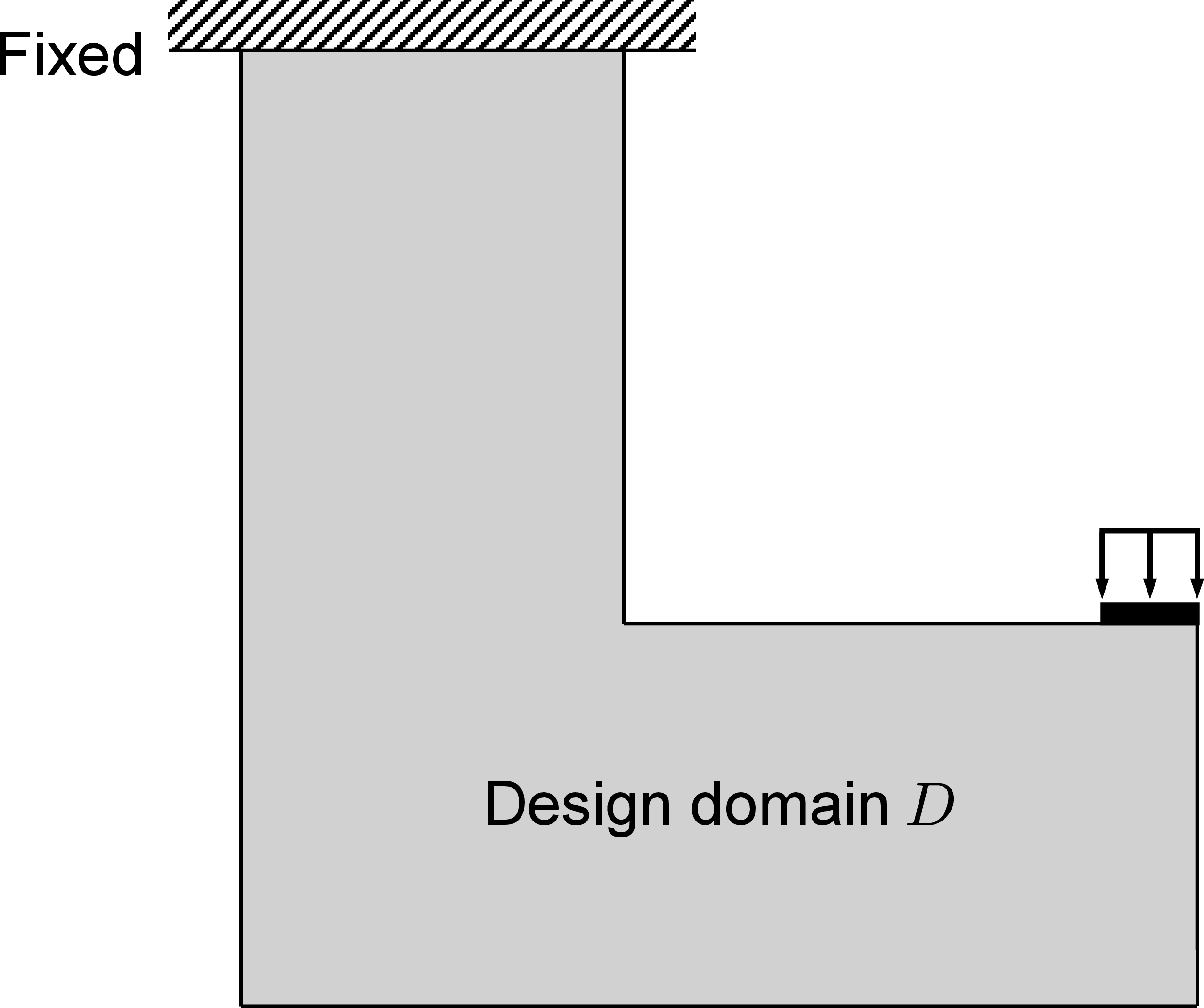}
    \caption{Design domain and boundary conditions of L-bracket}\label{fig_dd}
\end{figure}

As a numerical example, we solve the structural design problem of a two-dimensional L-bracket whose design domain and boundary conditions are shown in Fig.~\ref{fig_dd}.
It is widely used as a benchmark for stress-based optimization \citep{yang1996, duysinx1998, le2010, holmberg2013} and is known for causing strong nonlinearity due to stress concentration at the reentrant corner within the design domain.
Its optimization problem is formulated as a multi-objective optimization problem as follows:
\begin{eqnarray}\label{eq_str}
	\begin{aligned}
		& \underset{\boldsymbol{\rho}}{\text{minimize}}
			&& F_{1}=\underset{e}{\max}\left(\sigma_{e}\right),\\
			&&& F_{2}=\frac{\sum_{e=1}^{N}v_e\rho_e}{\sum_{e=1}^{N}v_e},\\
		& \text{subject to} 
			&& \rho_{e}\in\{0, 1\}\quad (e=1, 2, \ldots, {N}),
	\end{aligned}
\end{eqnarray}
where $\sigma_{e}$ is the von Mises stress in the $e$-th element.
$v_e$ and $N$ are the elemental volume and the number of elements, respectively.
In this paper, $N$ is set to 6400, indicating that the design domain is discretized into 6400 finite elements using structured meshes.
One of the objective $F_1$ is the maximum stress, making the optimization problem of \eqref{eq_str} a minimax one, while the other objective $F_2$ represents the volume fraction.
Note that each element $\rho_e$ of the design variable vector $\boldsymbol{\rho}$ takes discrete values of 0 or 1, representing a material distribution without intermediate densities known as grayscale.

While the optimization problem in \eqref{eq_str} is the original one to be solved, the low-fidelity optimization problem in \eqref{eq_strlf}, which is solved for initial data preparation and mutation in DDTD, is formulated as follows:
\begin{eqnarray}\label{eq_strlf}
	\begin{aligned}
		& \underset{\boldsymbol{\rho}^{(k)}}{\text{minimize}}
			&& \widetilde{F}=\boldsymbol{f}^{\mathrm{T}}\boldsymbol{u},\\
		& \text{subject to} 
			&& \boldsymbol{K}\boldsymbol{u}=\boldsymbol{f},\\
            &&& \widetilde{G}=\frac{\sum_{e=1}^{N}v_e\rho^{(k)}_e}{\sum_{e=1}^{N}v_e}-V_{\text{f}}^{\text{max}}\leq 0,\\
			&&& \rho^{(k)}_{e}\in[0, 1]\quad (e=1, 2, \ldots, {N}),\\
		& \text{for given}
			&& \boldsymbol{s}^{(k)},
	\end{aligned}
\end{eqnarray}
where the vectors $\boldsymbol{f}$ and $\boldsymbol{u}$ represent the external force and displacement, respectively, which form the equilibrium equations with the global stiffness matrix $\boldsymbol{K}$.
The objective $\widetilde{F}$ and constraint function $\widetilde{G}$ are the mean compliance and volume fraction, respectively.
The low-fidelity optimization problem of \eqref{eq_strlf} is a general stiffness maximization problem, which is easily solved using the density-based method \citep{bendsoe2003} with design variables $\rho^{(k)}_{e}$ relaxed to continuous values between 0 and 1.
In this paper, filter radius $r$ in density filter \citep{bruns2001, bourdin2001} and constraint values of volume fraction $V_{\text{f}}^{\text{max}}$ are employed as seeding parameters $\boldsymbol{s}$, i.e., they can be denoted as $\boldsymbol{s}=[r, V_{\text{f}}^{\text{max}}]$.

Table~\ref{tab_param} and \ref{tab_vae} list the parameters regarding the overall procedures of DDTD and the VAE, respectively.
Previous studies \citep{kii2024, kato2024} have demonstrated the effectiveness of DDTD in stress-based topology optimization, and this study investigates the impact of the proposed selection strategy on the solution search performance of DDTD under different parameter settings.
Among the various parameters, population size is known to significantly influence the search performance of GAs in the literature \citep{goldberg, koumousis2006}, thus we compare the optimization results with the three values shown in Table~\ref{tab_param}.
The prior study on DDTD \citep{kii2024} has also demonstrated that large population size leads to superior optimized solutions.

\begin{table}[t]
    \caption{Parameters for the overall procedures of DDTD}\label{tab_param}
    \begin{tabular*}{\textwidth}{@{\extracolsep\fill}lll}
        \toprule
        Description                                         & Symbol                        & Value\\
        \midrule
        Maximum iterations                                  & $t_{\text{max}}$              & 200\\
        Number of initial data                              & $N_{\text{ini}}$              & 100\\
        Population size                                     & $N_{\text{pop}}$              & 50, 100, 200\\
        Number of generated data in crossover using a VAE   & $N_{\text{VAE}}$              & 50, 100, 200 (aligned with $N_{\text{pop}}$)\\
        Number of mutants                                   & $N_{\text{mut}}$              & 16\\
        Iteration interval of mutation                      & $t_{\text{mut}}$              & 5\\
        Overlapping parameter for mutation                  & $G_{\text{mut}}^{\text{max}}$ & 0.01\\
        \botrule
    \end{tabular*}
\end{table}

\begin{table}[t]
    \caption{Parameters for the VAE and latent crossover}\label{tab_vae}
    \begin{tabular*}{\textwidth}{@{\extracolsep\fill}ll}
        \toprule
        Description                                 & Value\\
        \midrule
        Size of input and output Data               & 6400\\
        Size of latent space                        & 8\\
        Number of neurons for hidden layers         & 512\\
        Structure of the encoder network            & [6400, 512, 8]\\
        Structure of the decoder network            & [8, 512, 6400]\\ 
        Activation function for each layer          & Relu (hidden layers)\\
                                                    & Sigmoid (output layer)\\
        Optimizer                                   & Adam\\
        Reconstruction loss function                & Mean squared error\\
        Weight for Kullback-Leibler (KL) divergence & 0.001\\
        Number of epochs                            & 500\\
        Batch size                                  & 10\\
        Learning rate                               & 0.001\\
        Operator for latent crossover               & Simplex crossover (SPX) \citep{tsutsui1999}\\
        Number of parent individuals for SPX        & 9\\
        Expansion rate of simplex for SPX           & $\sqrt{10}$\\
        \botrule
    \end{tabular*}
\end{table}

\subsection{Verification of topological analysis using persistent homology}\label{sec42}
\begin{figure}[t]
    \centering
    \includegraphics[width=\textwidth]{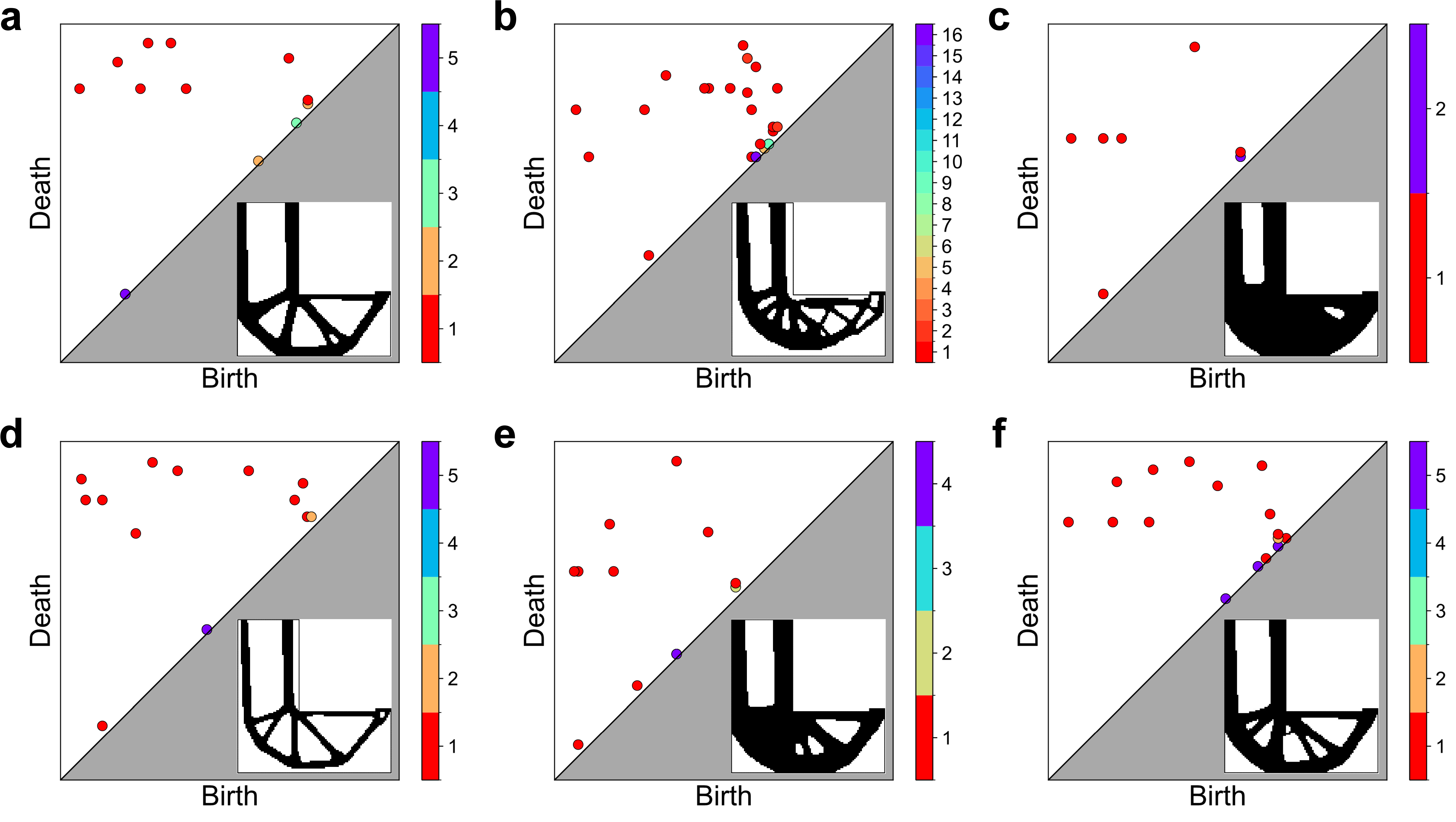}
    \caption{Examples of a pair of material distribution and its persistence diagram}\label{fig_result_pd}
\end{figure}

First, we verify whether the topological features are correctly extracted from the material distribution data using PH.
The material distributions of the initial data and mutants for them are binarized into black and white images with a threshold of $\rho^{(k)}_{e}=0.5$.
PH for these images was computed, and a part of the resulting PDs is shown in Fig.~\ref{fig_result_pd}.
Python software HomCloud \citep{obayashi2022} (version 4.4.1) was employed for the calculation.
The material distribution in Fig.~\ref{fig_result_pd}a confirms that the structure has a total of seven holes, including those consisting of the boundaries of the design domain.
Its PD shows seven points corresponding to these holes in the region far from the diagonal.
Similarly, for the more complex material distribution in Fig.~\ref{fig_result_pd}b with a greater number of holes, corresponding points can be observed on the PD.
All PDs in Fig.~\ref{fig_result_pd} show several plots near the diagonal, which can be regarded as noise and have little effect on the Wasserstein distance as described in Section~\ref{sec312}.
The PD calculated from the material distribution in Fig.~\ref{fig_result_pd}f, which includes two small holes added at the connections between each component in the material distribution of Fig.~\ref{fig_result_pd}a, shows two additional points corresponding to them.
Based on the definition of scale parameter shown in Fig.~\ref{fig_pds}, points further to the right on the PD correspond to smaller holes that appear earlier, and four additional plots in Fig.~\ref{fig_result_pd}d correspond to the small holes in the structure. 
The PDs corresponding to various structures, including Fig.~\ref{fig_result_pd}b with many holes and Fig.~\ref{fig_result_pd} with few holes, demonstrate that the number and size of the holes can be effectively captured, confirming that the topological features are analyzed using PH.

Next, we validate whether the Wasserstein distance between PDs is accurately measured.
Table~\ref{tab_w} compares the Wasserstein distance between PDs to the $L^2$ norm of material distributions shown in Fig.~\ref{fig_result_pd}.
Here we discuss the 2-Wasserstein distance calculated with $p=2$ in \eqref{eq_w}.
Based on the Wasserstein distances between PDs, the furthest material distributions are those in Fig.~\ref{fig_result_pd}a and \ref{fig_result_pd}c, as well as Fig.~\ref{fig_result_pd}c and \ref{fig_result_pd}d, while the closest material distributions are Fig.~\ref{fig_result_pd}a and \ref{fig_result_pd}d.
The $L^2$ norm similarly measures Fig.~\ref{fig_result_pd}c and \ref{fig_result_pd}d as a distant pair, while Fig.~\ref{fig_result_pd}a and \ref{fig_result_pd}b also have large values. 
On the other hand, The closest pairs are Fig.~\ref{fig_result_pd}a and \ref{fig_result_pd}f, as well as Fig.~\ref{fig_result_pd}c and \ref{fig_result_pd}e, resulting in completely different results compared to the Wasserstein distance.
In particular, focusing on Fig.~\ref{fig_result_pd}a and \ref{fig_result_pd}d, where their topology is similar but the position of each component is off by a few pixels, as shown in Fig.~\ref{fig_lp}, the Wasserstein distance assesses them as close, whereas the $L^2$ norm shows a large value, indicating that it does not accurately measure their similarity.
These results illustrate that the Wasserstein distance between PDs can appropriately measure the topological differences between material distributions.

\begin{table}[t]
    \caption{Computational results of 2-Wasserstein distance between persistence diagrams (PDs) and $L^2$ norm between material distributions shown in Fig.~\ref{fig_result_pd}}\label{tab_w}
    \begin{tabular*}{.6\textwidth}{@{\extracolsep\fill}lll}
        \toprule
        Pair of PDs & 2-Wasserstein distance & $L^2$ norm\\
        \midrule
        (a) \& (b)  & 46.05             & \textbf{51.72}\\
        (a) \& (c)  & \textbf{52.52}    & 42.97\\
        (a) \& (d)  & \textbf{20.14}    & 34.53\\
        (a) \& (e)  & 39.54             & 33.84\\
        (a) \& (f)  & 27.92             & \textbf{22.97}\\
        (b) \& (c)  & 41.46             & 40.60\\
        (b) \& (d)  & 45.99             & 39.95\\
        (b) \& (e)  & 36.52             & 40.93\\
        (b) \& (f)  & 33.67             & 47.70\\
        (c) \& (d)  & \textbf{53.41}    & \textbf{50.82}\\
        (c) \& (e)  & 28.93             & \textbf{21.45}\\
        (c) \& (f)  & 44.86             & 36.17\\
        (d) \& (e)  & 42.25             & 42.25\\
        (d) \& (f)  & 25.69             & 36.57\\
        (e) \& (f)  & 32.10             & 24.83\\
        \botrule
    \end{tabular*}
    \footnotetext{Note: Larger and smaller distances are shown in bold.}
\end{table}

\subsection{Validation of effectiveness of proposed selection strategy}\label{sec43}
Based on the verification of PH, we validate the effectiveness of the proposed selection strategy on the solution search performance of DDTD through comparing it with the conventional selection operation with the crowding distance sorting.
We use the hypervolume indicator \citep{shang2021}, which is a measure of diversity and convergence performance in multi-objective optimization, as a search performance metric of DDTD.
In the case of the two-objective optimization problem of \eqref{eq_str}, the hypervolume is calculated as the area formed by the non-dominated solutions of rank 1 and a predetermined reference point in the objective space.
Thus, a larger hypervolume value indicates a more advanced Pareto front.
Since the training process of VAEs involves randomness in DDTD, we compare the optimization results over ten trials in the three different population sizes shown in Table~\ref{tab_param}.

\begin{figure}[t]
    \centering
    \includegraphics[width=.6\textwidth]{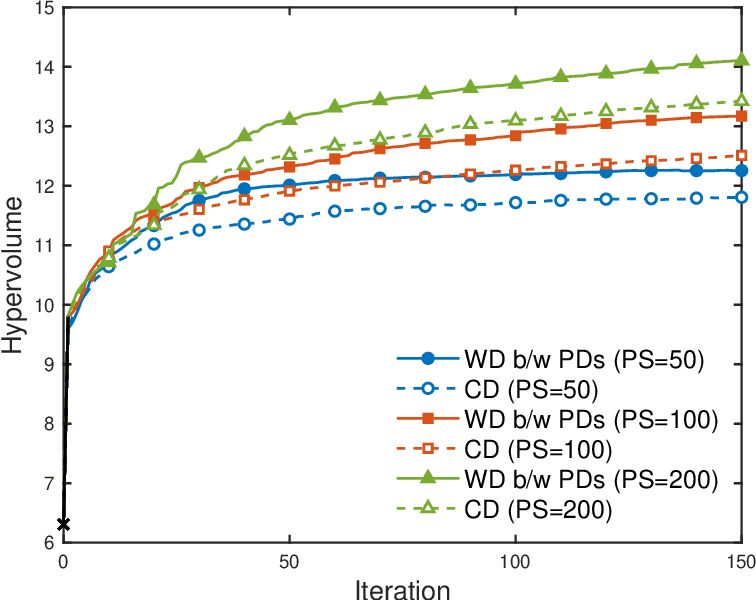}
    \caption{Iteration histories of the mean hypervolume over ten trials in different population size (PS):$N_{\text{pop}}=50, 100, 200$}\label{fig_hyv_itr}
\end{figure}

Figure~\ref{fig_hyv_itr} shows the iteration history of the mean hypervolume over ten trials.
Note that the reference point used for hypervolume calculations is common regardless of the population size.
Additionally, until iteration 1, the hypervolume values are nearly identical due to mutation, as represented by the black solid line in Fig.~\ref{fig_hyv_itr}.
In all cases, it can be confirmed that the proposed selection operation with the Wasserstein distance sorting between PDs outperforms the conventional one with the crowding distance sorting.
Quantitatively, the proposed method shows an increase of 8.19\%, 10.67\%, and 9.66\% over the conventional one for $N_{\text{pop}}=50, 100,$ and $200$, respectively.

\begin{figure}[t]
    \centering
    \includegraphics[width=.6\textwidth]{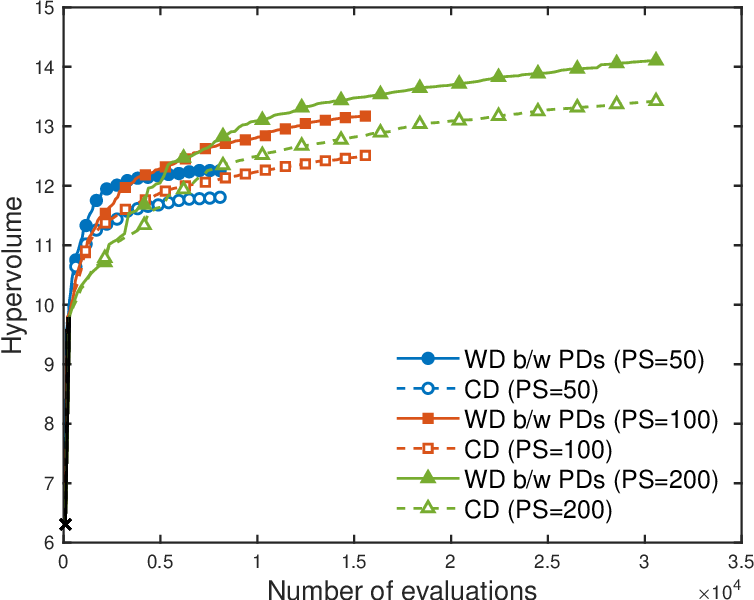}
    \caption{History of the mean hypervolume for the number of evaluations over ten trials in different population size (PS): $N_{\text{pop}}=50, 100, 200$}\label{fig_hyv_evl}
\end{figure}

Figure~\ref{fig_hyv_evl} shows the hypervolume history based on the number of evaluations.
Here, the number of evaluations refers to the count of performance evaluations of candidate solutions, which increments $N_{\text{pop}}+N_{\text{VAE}}$ by per iteration.
Although it is expected that the hypervolume value increases with a larger population size based on its definition, the final value for the proposed method with $N_{\text{pop}}=50$ is equivalent to that of the conventional method with $N_{\text{pop}}=100$.
Similarly, the proposed method with $N_{\text{pop}}=100$ requires only half the number of evaluations to yield comparable results comparable with the conventional method with $N_{\text{pop}}=200$.
The original paper \citep{yaji2022} also states that the most computationally expensive part of DDTD is the performance evaluation using the finite element method, and these results indicate that the proposed selection method can significantly reduce the computational time of DDTD.
It should be noted that the results suggest the potential usefulness of the proposed selection strategy for more complex topology optimization problems that involve higher computational costs for finite element analysis, such as three-dimensional problems or turbulent flow problems.

\begin{figure}[t]
    \centering
    \includegraphics[width=.8\textwidth]{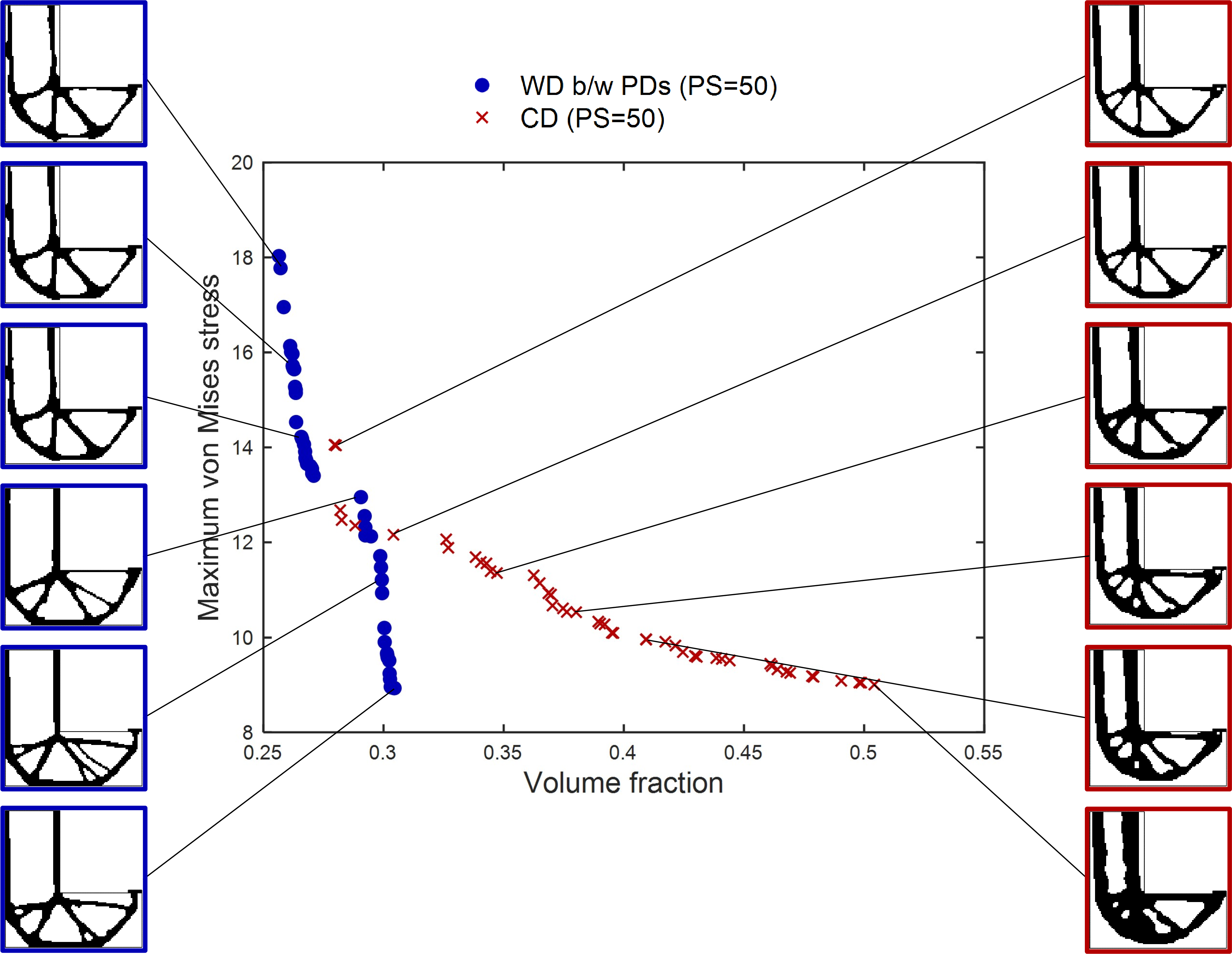}
    \caption{Objective space and some material distributions in optimization results for population size $N_{\text{pop}}=50$}\label{fig_pareto}
\end{figure}

\begin{figure}[t]
    \centering
    \includegraphics[width=\textwidth]{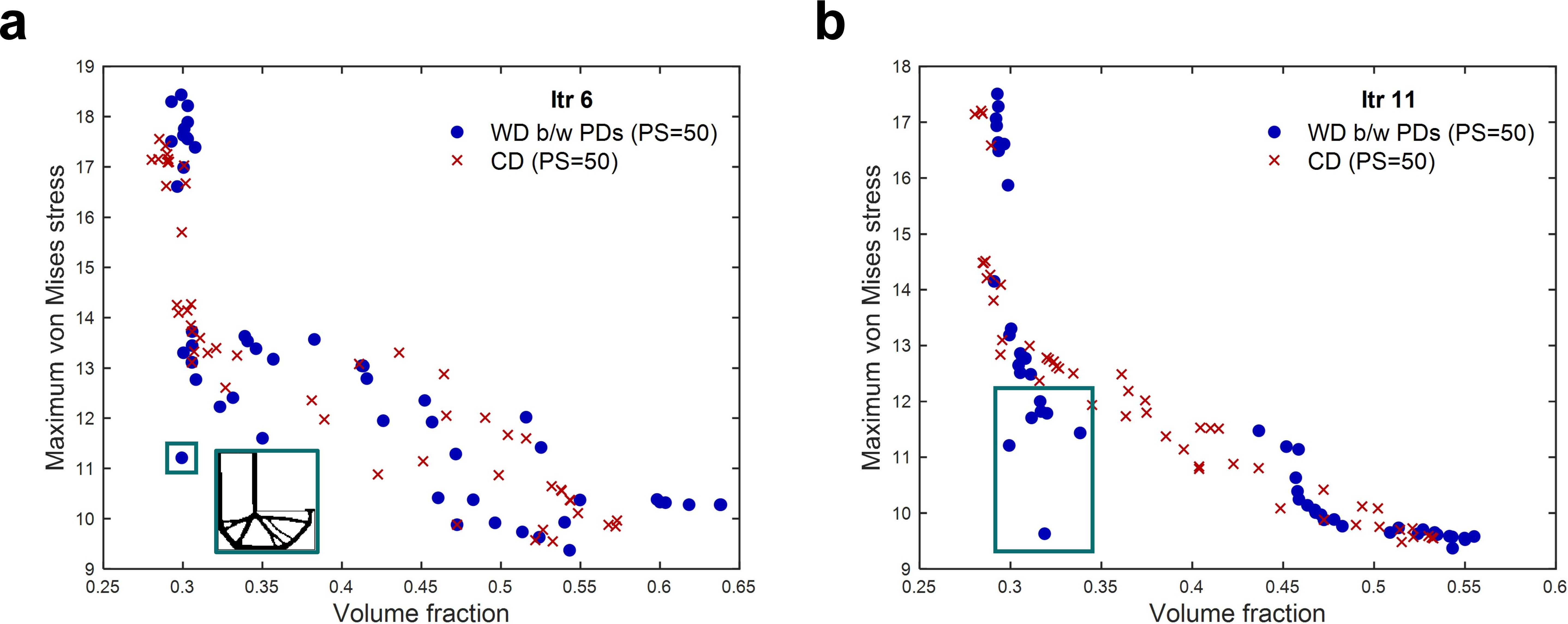}
    \caption{Objective space and extraction of superior solutions in the early stages of optimization (iterations 6 and 11) for population size $N_{\text{pop}}=50$}\label{fig_pareto_ini}
\end{figure}

For the case of population size $N_{\text{pop}}=50$, Fig.~\ref{fig_pareto} shows the Pareto front and some optimized structures from the trial with the maximum hypervolume value out of the ten trials illustrated in Figures~\ref{fig_hyv_itr} and \ref{fig_hyv_evl}.
As indicated by the hypervolume shown in Fig.~\ref{fig_hyv_itr}, the Pareto front obtained by the proposed method is much more advanced, especially in terms of volume reduction.
Focusing on material distributions, solutions with relatively large maximum stress values tend to have similar structures.
In contrast, solutions obtained by the proposed method with the maximum von Mises stress of less than 12 have unique structures, which contribute to the large hypervolume values.
Investigating when these unique solutions first appeared in the optimization calculations using the proposed method, they appeared as a mutant in iteration 6, as shown in Fig.~\ref{fig_pareto_ini}a.
Through subsequent generation, at iteration 11, it has multiplied through crossover, and an even superior solution with the maximum stress value of less than 10 has appeared, as shown in Fig.~\ref{fig_pareto_ini}b.
These results suggest that the proposed selection operation effectively enhances the population diversity, maintaining the population with a variety of design variables, and allowing crossover and mutation to produce novel superior ones that could not be achieved by conventional methods.
Additionally, as shown in Fig.~\ref{fig_pareto_ini}, the proposed Wasserstein distance sorting, unlike the crowding distance sorting, does not consider the proximity of solutions in the objective space during the early stage of optimization, resulting in scattered solution distributions with gaps in the front.
On the other hand, a continuous and uninterrupted Pareto front is eventually obtained as shown in Fig.~\ref{fig_pareto}, indicating that the proposed strategy of switching to the crowding distance sorting works correctly based on the theory of exploration and exploitation in EAs.

\section{Conclusion}\label{sec5}
This paper proposed a selection strategy enhancing the population diversity of solutions for data-driven topology design (DDTD).
Motivated by the need to consider the inherent diversity of material distributions in optimization problems with significant nonlinearity, we focused on persistent homology (PH) as a method for analyzing topology.
As a specific selection operation, we introduced the Wasserstein distance sorting between persistence diagrams instead of the crowding distance sorting in the non-dominated sorting genetic algorithm II (NSGA-II), a type of evolutionary algorithm for multi-objective optimization problems.
In the numerical example of stress-based topology optimization, it was confirmed that PH effectively analyzes the holes in material distribution data and that the Wasserstein distance between persistence diagrams is appropriately calculated.
It was demonstrated that the proposed selection operation improves the solution search performance of DDTD and leads to the discovery of unique and high-performance structures.

One of the significant achievements of this paper is demonstrating that solution search performance in DDTD is not compromised even with reduced population size.
Our future work will focus on tackling large-scale or complex optimization problems that involve high computational costs for physical performance analysis.

\backmatter





\bmhead{Acknowledgement}
This work was supported by JSPS KAKENHI Grant Numbers 23H03799 and 24KJ1640.

\section*{Declarations}
\bmhead{Conflict of interest}
The authors declare that they have no conflict of interest.

\bmhead{Replication of results}
The necessary information for replication of the results are presented in the manuscript.
Interested readers may contact the corresponding author for further details regarding the implementation.

\bibliography{references}

\end{document}